\documentclass{amsart}
\usepackage{amssymb}
\newtheorem{theorem}{Theorem}[section]
\newtheorem{lemma}[theorem]{Lemma}

\newtheorem{corollary}[theorem]{Corollary}
\newtheorem{example}[theorem]{Example}

\newtheorem{remarks}[theorem]{Remarks}
\newtheorem{remark}[theorem]{Remark}
\newtheorem{definition}[theorem]{Definition}
\newcommand{\calT}{\mbox{${\mathcal T.}$}}
\newcommand{\calt}{\mbox{${\mathcal T}$}}
\newcommand{\calu}{\mbox{${\mathcal U}$}}

\newcommand{\calW}{\mbox{${\mathcal W}$}}

\newcommand{\calB}{\mbox{${\mathcal B}$}}
\newcommand{\calU}{\mbox{${\mathcal U}$}}
\newcommand{\calV}{\mbox{${\mathcal V}$}}

\newcommand{\BBB}{\mathbb B}
\newcommand{\BBN}{\mathbb N}
\newcommand{\BBQ}{\mathbb Q}
\newcommand{\BBP}{\mathbb P}
\newcommand{\BBR}{\mathbb R}
\newcommand{\BBZ}{\mathbb Z}
\begin{document}
\begin{center}{\Large Looking through newly to the amazing irrationals}\\\smallskip
{\normalsize Pratip Chakraborty\footnote{\copyright\,My pleasure if you need it}\\
University of Kaiserslautern\\Kaiserslautern, DE
67663\\chakrabo@mathematik.uni-kl.de
 \\ 14/12/2004}

\end{center}\bigskip
\section{Here we start with metric spaces and dimension zero}
We assume basic familiarity with the theory of metric spaces, in
particular, the space of real numbers with its usual metric. In this
section, we review some notation and  a few basic definitions and
theorems concerning metric spaces (including the concept of
dimension zero).  In \S 3 we use the completeness (or least upper
bound) axiom of the real numbers in the form that says the space of
real numbers is locally compact.

\subsection{Basic definitions}

We use standard notation: $\BBR$ denotes the set of real numbers,
$\BBP$ the set of irrational numbers, $\BBQ$ the set of rational
numbers, $\BBZ$ the  set of all integers, $\BBN$ the set of all
positive integers, and $\omega$ the set of non-negative integers
(the first infinite ordinal).
\begin{definition}\label{metricdef}A pair $(X,d)$ is called a {\em metric space} provided $X$ is a
set and $d$ is a {\em metric}, i.e., a function $d:X\times
X\rightarrow [0,\infty)$ which satisfies the following three
properties for every x,y,z in X:

\smallskip
(a)  $d(x,y) = 0$ if and only if $x = y$

\smallskip
(b)  $d(x,y) = d(y,x)$

\smallskip
(c)  $d(x,z) \leq d(x,y) + d(y,z)$ {\em (triangle inequality)}

\end{definition}

Let $(X,d)$ be a metric space, $x\in X$, and $r > 0$.  The set
$$N(x,r) = \{y\in X: d(x,y) < r\}$$
is called a {\em basic neighborhood} of $x $ (of radius $r$).  A set
$U\subset X$ is called an {\em open set} provided $U$ is a union of
basic neighborhoods, i.e., for every $x\in U$ there exists $r = r(x)
> 0$ such that $N(x,r)\subset U$.

The set

\[T = T(X,d) = \{U\subset X: U \mbox{ is an open set}\}\]
is called the {\em topology of $(X,d)$}, or the {\em topology on $X$
induced by the metric~$d$}. A family $\calB\subset \calT$ is called
a {\em base} for the topology $\calt$ provided every element of
$\calT$ is a union of elements from $\calB$.   A set $F\subset X$ is
called {\em closed} (in $X$) provided $X-F$ is open in $X$.  If $U$
is open and $F$ is closed, then $U\setminus F$ is open and
$F\setminus U$ is closed (exercise).

For example, the function $d(x,y) = |x-y|$ (the absolute value) can
easily be seen to be a metric on $\BBR$, and is called the {\em
usual metric} on $\BBR$.

If $(X,d)$ is a metric space and $Y \subset X$, the {\em induced
metric} on $Y$ is $d|(Y\times Y)$.  Thus $(Y,d|(Y\times Y))$ is a
metric space  and is called a {\em subspace} of $(X,d)$.  For
example, the set of rational numbers $\BBQ$ is a subspace of the
metric space  $\BBR$, where
  $\BBR$ has its usual metric, and
likewise the set of irrational numbers  $\BBP$, is a subspace of
$\BBR$.  In both
 $\BBQ$ and  $\BBP$
there are many sets which are both open and closed at the same time,
and these sets play an important role in these spaces. The
properties of a set being open or closed are relative to subspaces,
and their relation to subspaces is a simple one:  Let $(X,d)$ be a
metric space and $Y \subset X$.  If $y \in Y$ and $r > 0$, let
$N_Y(y,r)$ denote the basic neighborhood of $y$ in the space $Y$,
and let $N_X(y,r)$ denote the basic neighborhood of $y$ in $X$.
Then $N_Y(y,r) = N_X(y,r)\cap Y$.  It follows that a set $U \subset
Y$
 is open in $Y$ if
and only if there exists an open set $V$ in $X$ such that $U = V\cap
Y$.  An analogous statement holds for closed sets. For short, we say
that a set $A \subset X$ is {\em clopen} provided $A$ is both open
and closed in $X$.  If $U \subset Y\subset X$, then $U$
 can be open (or
closed) in $Y$ without being open (or closed ) in $X$.
 For example if $a,b \in \BBQ$ with $a < b$  then $(a,b)\cap \BBP$ is
both open and closed in  $\BBP$, but neither open nor closed in
$\BBR$.

A family $\calB$ of open subsets of $X$ is called a {\em base for
$X$} if every open set is a union of elements from  $\calB$.  The
set of all open intervals with rational end-points is a countable
base of $\BBR$, and when intersected with $\BBP$ forms a countable
base for $\BBP$. A set $D\subset X$ is called {\em dense} provided
every non-empty open set contains a member of $D$. A metric space
$X$ is called {\em separable} if $X$ has a countable dense set. It
is well-known that a metric space is separable if and only if it has
a countable base (exercise). We speak of a separable metric space
even when we have ``a countable base'' in mind.

\subsection{Homeomorphisms}
The study of the irrational numbers involves consideration of spaces
which are ``identical in the sense of topology" to the space of
irrational numbers (for instance the Baire space considered in the
next section). We recall the basic concepts of
 ``homeomorphism," which is the precise notion
of ``identical in the sense of topology."  We also mention a
special kind of homeomorphism called topological equivalence.

Let $(X,d)$ and $(Y,\rho)$ be metric spaces, and $f:X\rightarrow Y$
 a function.  We say $f$ is {\em continuous} provided
$f^{-1}(V) = \{x \in X: f(x) \in V\}$ is open in $X$ for every open
set $V \subset Y$ (continuity can be stated, of course, in terms of
the two metrics $d$ and $\rho$).

A function that is both one-to-one and onto is called a {\em
bijection}.  A bijection $f$ such that $f$ and $f^{-1}$ are
continuous is called {\em bicontinuous}.

\begin{definition} A function $h:X\rightarrow Y$ is called a {\em homeomorphism} (between the two
metric spaces $(X,d)$ and $(Y,\rho)$) provided (i) $h$ is a
bijection, (ii) $ h$ is continuous,
 and (iii) $h^{-1}$ is continuous (i.e., $h$ is bicontinuous).  We also say that $X$ and $Y$ are {\em homeomorphic},
or $X$ {\em is homeomorphic  to} $Y$.
\end{definition}

\begin{example}\label{twointervals}  Let $X = (a,b)$ and $Y = (c,d)$ be two open intervals in
$\BBR$ with the usual metric, then $X$ and $Y$ are homeomorphic.
\end{example}
Proof.  Let $h$ be the linear map whose graph is a straight line in
the Euclidean plane through the two points in the plane $\langle
a,d\rangle$ and $\langle b,c\rangle$. Then $h$ is a bijection, and
is continuous (even differentiable), and $h^{-1}$ is continuous
since its graph is also a straight line in the plane.

\begin{example}\label{intervals}.  The real line $\BBR$ is homeomorphic to the open
interval $(-1,1) $ (and therefore $\BBR$ is homeomorphic to every
open interval by \ref{twointervals}).
\end{example}
Proof.  Define $h:\BBR\rightarrow (-1,1)$ by $h(x) = x/(1+|x|)$.  To
see that $h^{-1}$ is continuous note that $h^{-1}: (-1,1)
\rightarrow \BBR$ is defined by   $h^{-1}(x) = x/(1-|x|)$. For an
alternate proof that $\BBR$ is homeomorphic to every open interval,
use the function $\arctan(x)$ which maps $\BBR$ onto $(-\pi/2,
\pi/2)$.

 \medskip
 A property is called {\em invariant under homeomorphisms} provided for any two
homeomorphic spaces, one has the property if and only if the other
has the property. A property which is invariant under homeomorphisms
is called a {\em topological property}.
 Loosely speaking, a property
which can be stated in terms of open sets, without mentioning the
metric, is a topological property (two such properties are
compactness---not discussed here--- and separability).  Any property
which is preserved by continuity, is a topological property.  An
example of a property which is not a topological
 property is ``boundedness." For example,  $(-1,1)$ and $\BBR$ are homeomorphic, but only
one is bounded.

\begin{definition}If X is a set and $d, \rho$ are two metrics on $X$, we say that $d$ and $\rho$ are
{\em topologically equivalent} provided $T(X,d) = T(X,\rho)$, in
other words, provided the identity function on $X,
id:(X,d)\rightarrow (X,\rho)$ is a homeomorphism.
\end{definition}
We use several times the obvious fact that the composition of two
homeomorphisms is a homeomorphism.

\subsection{Two Notions of Dimension Zero}
We conclude this section with a brief discussion of two notions of
dimension zero in metric spaces.  First we recall some terminology.

\begin{definition}  Let $(X,d)$ be a metric space, and $\calU, \calV$ two families of
subsets of $X$. We say that $\calV $ {\em refines} $\calU$ (or that
$\calV$ is a {\em refinement} of $\calU$) provided for every $V \in
\calV$
 there  exists $U \in \calU$ such that $V \subset U$.  A family \calU\
 is called {\em pairwise disjoint}
provided that for every $U,U' \in\calU$, if $U \not= U^{\prime}$,
then $U\cap U' = \emptyset$.  We say that $\calV$ is an {\em open
refinement of} $\calu$ provided $\calV$ is a family of open sets
that refines $\calu$, and covers the same set that $\calu$ covers
($\cup\calV =\cup\calu$).
\end{definition}

    The family $\calU = \{(n, n+1)\cap \BBP: n \mbox{ is an integer}\}$, is a
example of a pairwise disjoint clopen cover of the irrational
numbers, and $\calV = \{(n/2, (n+1)/2)\cap \BBP: n \mbox{ is an
integer}\}$ is a pairwise disjoint clopen cover of $\BBP$ such that
\calV\ refines \calU.

\begin{definition}\label{zerodim} A metric space $(X,d)$ is called {\em 0-dimensional} provided
for every $x \in X $ and every $r > 0$, there exists a set $U$ which
is clopen and $x \in U \subset N(x,r)$. A metric space $(X,d)$ is
said to have {\em covering dimension zero} provided that for every
open cover $\calU$ of $X$, there is a pairwise disjoint open cover
$\calV$ of $X$ such that $\calV$ refines \calU.
\end{definition}

    It is easy to see that every space with covering dimension zero
is 0-dimensional  (exercise), but the converse is not true. There
is a famous example of a metric space which is 0-dimensional
 but does not have covering
dimension zero, given in 1962 by Prabir Roy  \cite{roy}.  Other such
examples are now known.  Two easier such example, were given in 1990
by John Kulesza \cite{kulesza} and Adam Ostaszewski \cite{adam}.

It is often easier to see if a space is 0-dimensional, than if it
has covering dimension zero. For example, both the rational and
irrational numbers, with their usual metrics, are obviously
0-dimensional. We now give a result that implies  $\BBP$ and $\BBQ$
 have covering dimension zero.

\begin{lemma}\label{clopenbase} Every 0-dimensional separable metric space $X$
has a countable base consisting of clopen sets.  In particular,
$\BBP$ has such a base.
\end{lemma}

 Proof.  For every $x\in X$, and every $\epsilon > 0$, there exists a clopen set
$B(x,\epsilon)$ such that $x\in B(x,\epsilon)\subset N(x,\epsilon)$.
Clearly, $\{B(x,\epsilon): x \in X,\epsilon > 0\}$ is a  base for
$X$ consisting of clopen sets.  Since $X$ has a countable base, we
know that {\em every} base contains a countable subset which is also
a base (exercise). Thus there is a countable base consisting of a
subset of these $B(x,\epsilon)$'s, and that completes the proof.

\begin{theorem}\label{0covering} If $(X,d)$ is  a separable metric space, then $(X,d)$ is 0-dimensional
if and only if $(X,d)$ has covering dimension zero.  In particular,
$\BBP$ has covering dimension zero.
\end{theorem}

    Proof. Given the exercise following Definition \ref{zerodim}, we need only show
that if  $X$ is 0-dimensional and separable, then $X$ has covering
dimension zero.
 By \ref{clopenbase}, $X$ has a countable base $\calB$ consisting of clopen sets.
   Let $\calU$ be an open cover
of $X$, and let
$$\calW = \{B \in \calB: B \subset  U \mbox{ for some } U \in \calU\}.$$
Since $\calW$ is countable, put $\calW = \{B_n: n \in \omega\}$.
Now $\calW$ is an open refinement of $\calU$, but may not be
pairwise disjoint. We define a sequence of open sets by induction:
Put $V_0 = B_0$, and for $n \geq 1$ define $V_n =
B_n\setminus\cup\{B_i: i < n\}$.  Clearly $\calV = \{V_n: n \in
\omega\}$ is a family of pairwise disjoint open (and closed) sets,
each of which is contained in some member of $\calU$. To complete
the proof we need only show that $\calV$ covers $X$, and this is
clear since for every $x \in X$ there is a first $n$ such that $x
\in B_n$; so $x \in V_n$.  Thus $\calU$ has a pairwise disjoint open
refinement.

\medskip
If $\calU$ is a nonempty family of subsets of a metric space $X$, we
recall that the {\em diameter} of a set $U$ is defined by
$\mbox{diam\,}(U) =\sup\{d(x,y): x,y \in U\}$, and {\em mesh} of a
family of sets $\calU$ is defined by $\mbox{mesh\,}(\calu) =
\sup\{\mbox{diam\,}(U): U \in \calU\}.$

\section{Baire space and ultrametrics}

An important space which is homeomorphic to the irrationals is the
Baire space  $\mathbb B$. Since the natural metric on  $mathbb B$
is an ultrametric, we study these metrics in some detail.  In
particular, we give a topological characterization of ultrametrics
(\ref{ultrametricchar}).  Ultrametric spaces, like  $\mathbb B$,
have metric properties which are very different from the properties
of the usual metric on Euclidean (sub)spaces (e.g., $\BBR, \BBP$,
and $\BBQ$); see \ref{ultrametric}.

\subsection{The Baire space}
Recall that $\omega$ denotes the first infinite ordinal (i.e., the
set of non-negative integers), and let $^\omega\omega$ denote the
set of all functions $f:\omega\rightarrow \omega$.  We can think of
a point $f$ in $\BBB$ as a sequence of non-negative integers $f =
(f_0,f_1,\cdots , f_n,\cdots ).$

\begin{definition}\label{BaireMetric}For two distinct functions
$f,g  \in    ~^\omega\omega~$, let $k(f,g) $ denote the first $n
 \in    \omega$ such that $f(n) \not=  g(n)$.  We define a metric $d$ on
$^\omega\omega$ by
  \[d(f,g) = \left\{ \begin{array}{ll}
1/(k(f,g)+1) & \mbox{if $f\not= g$}\\
0 & \mbox{if $f = g$}
\end{array}\right.\]
\end{definition}
\begin{lemma}\label{Bairemetriclemma}  The function $d$ is a metric on
$^\omega\omega$.
\end{lemma}

Proof.  We need only check the triangle inequality since the other
properties are obvious.  Let $f,g,h  \in    ~^\omega\omega~$.  We
must show that $d(f,h) \leq   d(f,g) + d(g,h)$.  In fact we will
show that
$$d(f,h) \leq   \max\{d(f,g), d(g,h)\}.$$
If $f = h, f = g$ or $g = h$, this is trivial; so assume otherwise
and let  $n = k(f,h), m = k(f,g)$ and $q = k(g,h)$.  Thus, it
suffices to show that either $1/(n+1) \leq 1/(m+1)$, or $1/(n+1)
\leq   1/(q+1)$; or, equivalently, to show that either $m \leq   n$
or $q \leq   n$.  Suppose that $m > n$; then $f$ and $g$ agree past
$n = k(f,h)$; so $g(n) =f(n) \not=  h(n)$.  But this says that
$k(g,h) \leq   n$, i.e., $q \leq   n$.

\begin{definition}The metric space $ \BBB     = (^\omega\omega, d)$, where $d$ is defined
in \ref{BaireMetric}, is called the {\em Baire space}.
\end{definition}

The basic neighborhoods $N(f,\epsilon)$ in the Baire space can be
described in another way. For each $n  \in    \omega$, let
$^{n}\omega$ denote the set of all functions from $n$ into $\omega$,
and let $^{<\omega}\omega =  \cup  \{^{n}\omega: n \in    \omega\}$.
For each  $\sigma     \in    ^{n}\omega$, define $[ \sigma   ] = \{f
\in ^\omega\omega:  \sigma     \subset   f\}$. It follows from the
next lemma that the basic neighborhoods in  $\BBB\ $    are the same
sets as the $[ \sigma   ] $ for $\sigma \in ~^{<\omega}\omega$.
These sets $[\sigma]$, of course, form the natural base for the
product topology on $\BBB$.

 For $f  \in
^\omega\omega$ and $m  \in    \omega$, we use the notation $f|m =
\{(i,j)  \in    f: i < m\}$; thus $f|m \in ~^{m}\omega~$

\begin{lemma}  For each  $n$ let
$f  \in ~^{\omega}\omega$, and $r > 0$. If $r > 1$, then in $\BBB$,
$N(r,f) = ~^{\omega}\omega~$.  If $r \leq  1$, then $N(f,r) =
[f|m]$, where $m \geq  2$ is the first natural number such that $1/m
< r \leq  1/(m-1)$.
\end{lemma}

Proof.  The proof is left as an exercise.

We want to consider a different representation of Baire space. Let
$^{\BBN}\!\BBN$ denote the set of all functions
$f:\BBN\rightarrow\BBN.$  Define
$$\BBB_2 = (\BBZ\times ^{\BBN}\!\BBN) =\{<z,h>: z \in \BBZ, \mbox{ and }h\in ^{\BBN}\!\BBN\}$$
A point $p \in \BBB_2$ can be considered as a sequence
$(p_0,p_1,\cdots p_n,\cdots)$ where $p_0\in \BBZ$ and $p_n\in\BBN$
for $n \geq 1$. If $p \not=q$ in $\BBB_2$, let $k_2(p,q)$ denote the
first integer $n \geq 0$ such that $p_n \not= q_n$.

\begin{definition}\label{BaireMetric2}For
$p,q  \in    \BBB_2$, define
  \[d_2(p,q) = \left\{ \begin{array}{ll}
1/(k_2(p,q)+1) & \mbox{if $p\not= q$}\\
0 & \mbox{if $p = q$}
\end{array}\right.\]
\end{definition}
\begin{lemma}\label{Bairemetric2}  The function $d_2$ is a metric on
$\BBB_2$.
\end{lemma}
Proof. This is similar to  Lemma \ref{Bairemetriclemma}.

\begin{theorem}\label{prod}$\BBB$ and $\BBB_2$ are homeomorphic.
\end{theorem}

Proof.  Let $\phi:\omega\rightarrow\BBZ$ and $\psi:\omega
\rightarrow\BBN$
 be any bijections.  Using the notation $f= (f_0,f_1,\cdots ,f_n,\cdots)$ for
a point in $\BBB$, define a function $\Psi:\BBB\rightarrow\BBB_2$ by
$$\Psi(f) = (\phi(f_0),\psi(f_1),\cdots ,\psi(f_n),\cdots ).$$
Then $\Psi$ is one-one and onto (exercise).  Moreover, $\Psi$ can be
seen to be bicontinuous from the following observation. If $f,g\in
\BBB$, and $p = \Psi(f), q =\Psi(g)$ then the first place where
$f,g$ differ is the same as the first place that $p,q$ differ.  Thus
for $f,g$ distinct,
$$d_2(p,q) = d_2(\Psi(f),\Psi(g)) = \frac{1}{k_2(\Psi(f),\Psi(g))+1}=
\frac{1}{k(f,g)+1}= d(f,g).$$

\medskip
As an aside, we mention that two metric spaces $(X,d)$ and
$(Y,\rho)$ are said to be {\em isometric} provided there exists a
bijection $ h:X\rightarrow Y$ such that for all $a,b \in X$ we have
$d(a,b) = \rho(h(a),h(b))$.  Such a function $h$ is called an {\em
isometry}. Clearly every isometry is a homeomorphism, but not every
homeomorphism is an isometry.

\medskip
In the next section, we will show that $\BBB_2$ (hence $\BBB$) is
homeomorphic to the irrational numbers.

\subsection{Ultrametrics}

\begin{definition} A metric $d$ on a set $X$ is called an {\em ultrametric (or a non-
Archimedean metric)} on $X$ (and $(X,d)$ is called an {\em
ultrametric (or non-Archimedean) space)\/} provided $d$ satisfies
the following stronger form of the triangle inequality: for all
$x,y,z  \in    X$,

\[d(x,z) \leq  \max\{d(x,y), d(y,z)\}\mbox{  (strong triangle inequality)}.\]
\end{definition}
The strong triangle inequality might be called the ``isosceles
triangle inequality" because it implies that among the three sides
of a triangle $\triangle xyz$, at least two have equal length, where
by the lengths of the  sides of $\triangle xyz$, we mean the numbers
$d(x,y), d(y,z)$ and $d(z,x)$. Since each of these three lengths is
less than or equal to the maximum of the other two, it follows that
two of these lengths are the same.

The proof of 1.3, shows that the natural metric defined on
$^\omega\omega$ above is an ultrametric; so  $\BBB\ $    is an
ultrametric space.  Note that the usual metric on   $\BBP \ $ (the
subspace of irrationals), is not an ultrametric (e.g., take  $x
=\sqrt{2}, y = 1+\sqrt{2}$, and $z = 2+\sqrt{2})$.  The following
results show how very different an ultrametric is from the usual
Euclidean metrics.

\begin{theorem}\label{ultrametric}  Let $(X,d)$ be an ultrametric space.  Then the following
hold:

(a)  For every $x \in X$ and every $r > 0$, $N(x,r)$ is an open and
a closed set.

(b)  If two basic neighborhoods have a point in common, then one is
contained in the other (i.e., if $N(x,r) \cap N(y,s) \not=
\emptyset$ and    $r \leq s$, then $N(x,r)  \subset   N(y,s)$.

(c)  If $N(x,r) \cap N(y,r) \not= \emptyset$, then $N(x,r) =
N(y,r)$.

(d)  Every point in a basic neighborhood qualifies as its ``center."
I.e., if $y  \in    N(x,r)$, then $N(x,r) = N(y,r)$.

(e)  Let $S(y,r) = \{x  \in    X: d(x,y) \leq r\}$.  If $x  \in
S(y,r)$ then $N(x,r) \subset S(y,r)$.

(f)  Every non-empty intersection of finitely many basic
neighborhoods is a basic neighborhood.

(g)  Every open set is a union of a family of pairwise disjoint
basic neighborhoods.

(h)  Any union of basic neighborhoods of the same radius is clopen.
\end{theorem}

Proof.  For (a), $N(x,r)$ is open by definition of the metric
topology.  To see that $N(x,r)$ is closed, we show that $X\setminus
N(x,r)$ is open. Let $y  \in    X\setminus N(x,r)$; then $d(x,y)
\geq  r$.  We claim that $N(y,r) \cap N(x,r) = \emptyset$.  If this
is not true, then there exists a point $z  \in    N(y,r) \cap
N(x,r)$; so $d(z,y) < r$ and $d(z,x) < r$.  By the strong triangle
inequality we have

$$r \leq  d(x,y) \leq  \max\{d(x,z), d(z,y)\} < r,$$

but this is clearly a contradiction.  Thus $N(y,r) \cap N(x,r) =
\emptyset$, and it follows that $N(x,r)$ is closed.

The remaining parts of this Theorem are left as exercises.

\subsection{A characterization of ultrametrizable spaces}

We now come to the main result of this section.

\begin{theorem}\label{ultrametricchar}  Let $(X,d)$ be a metric space.  The following three
conditions are equivalent.

1.  There exists an ultrametric  $\rho$  on $X$ which is
topologically equivalent to $d$.

2.  $(X,d)$ has covering dimension zero.

3.  There exists a sequence $\{\calU_i: i  \in    \omega\}$ of
pairwise disjoint open covers of $X$ such that

(i) if $U  \in    \calU_i$, then $diam(U) \leq  2^{-(i+1)}$, and

(ii) $\calU_{i+1}$ refines $\calU_i$, for all $i  \in    \omega$.

\end{theorem}
    Proof.  Proof of $1\rightarrow 2$.  Let  $\rho $  be an ultrametric on
$X$ such that
 $\rho$  and $d$ are topologically equivalent.  Let $\calU$ be an open cover
of $X$ (since  $\rho$  and $d $ are topologically equivalent the
term ``open" is not ambiguous).  We need to find a cover ${\mathcal
V}$ of pairwise disjoint open sets such that ${\mathcal V}$ refines
${\mathcal U}$.  For each $x$ in $X$, let

$$i(x) = \min\{i \geq 1: \exists U  \in    \calU
\mbox{ such that }N_{\rho} (x,1/i)  \subset   U\}.$$

Put $\calV = \{N_{\rho} (x,1/i(x)): x  \in    X\}$.  Clearly \calV\
is an open cover of $X$ which refines $\calU$; so we need only show
that \calV\ is pairwise disjoint.  Let $N_{\rho}(x,1/i(x))$ and
$N_{\rho}(y,1/i(y))$ be distinct elements of \calV, and suppose that
they are not disjoint; say $z  \in    N_{\rho}(x,1/i(x))  \cap
N_{\rho}(y,1/i(y))$.  By 1.7(b), one of these basic neighborhoods is
contained in the other; say $N_{\rho}(x,1/i(x))  \subset
N_{\rho}(y,1/i(y))$.  By definition of $i(y)$ there exists $U  \in
\calU$ such that $N_{\rho}(y,1/i(y))  \subset   U$. Since $x$ is in
$N_{\rho}(y,1/i(y))$, it can be used as the ``center", i.e.,
$N_{\rho}(y,1/i(y)) = N_{\rho}(x,1/i(y))$.  Thus we have

$$N_{\rho}(x,1/i(y)) = N_{\rho} (y,1/i(y)) \subset   U.$$

By the definition of $i(x)$, we have $i(x) \leq  i(y)$; so

$$N_{\rho}(x,1/i(x)) \supset N_{\rho}(x,1/i(y)) = N_{\rho}(y,1/i(y))$$

Hence $N_{\rho}(x,1/i(x)) = N_{\rho}(y,1/i(y))$ which contradicts
our assumption that these two elements of $V$ were distinct.

Proof of $ 2\rightarrow 3$.  By hypothesis, $X$ has covering
dimension zero. Let $\calU_0$ be a pairwise disjoint open refinement
of $\{N(x,1) : x  \in X\}$. Assume that we have defined pairwise
disjoint open covers $\calu_i$, for $i \leq  n$ such that properties
(i) and (ii) in (3) hold for all $i \leq  n$.  We construct
${\mathcal U}_{n+1}$ as follows.  Let

$${\mathcal V}=\{N(x,2^{-(n+3)}): x  \in    X\}\mbox{, and }
{\mathcal W} = \{ N  \cap  U: N  \in {\mathcal V} \mbox{ and } U
\in    {\mathcal U}_n\}.$$

Then ${\mathcal W}$ is an open cover of $X$; so by hypothesis, there
exists a pairwise disjoint open refinement of $W$, call it
${\mathcal U}_{n+1}$.  It is easy to see that the sequence
$\{\calU_i: i  \in    \omega\}$ satisfies the requirements of (3).

Proof of $3\rightarrow 1$.  Let $\{\calU_i: i  \in    \omega\}$ be a
sequence of open covers satisfying the properties in (3). We define
an ultrametric  $\rho$  on $X$ as follows:  For any two distinct
points $x, y \in X$ define

$$k(x,y) = \min\{i  \in    \omega: x \mbox{ and } y \mbox{ are in different
members of } {\mathcal U}_i\}.$$

Property (i) implies that $k(x,y)$ is well-defined.  For any $x, y
\in X$ define (and compare with \ref{BaireMetric})

  \[\rho(x,y) = \left\{ \begin{array}{ll}
1/(k(x,y)+1) & \mbox{if $x\not= y$}\\
0 & \mbox{if $x = y$}
\end{array}\right.\]
We have to check first that  $\rho$  is an ultrametric, and the only
part of the definition that is not obvious is the stronger form of
the triangle inequality.  Let, $x,y,z$ be in $X$.  We will show that
$$\rho (x,y) \leq  \max\{ \rho (x,z),  \rho (z,y)\}.$$
If $x = y$ this is trivial; so we assume that $x \not=  y$. By (i)
there exists a first $n \geq 0$ such that $x$ and $y$ are in
different members of $\calU_n$; so $k(x,y) = n$.  Let $U, U'  \in
\calU_n$ such that $x  \in U$ and $y  \in    U'$. If $z$ is not in
$U$; then $z$ and $x $ are in different members of $\calU_n$ which
implies that $k(z,x) \leq n$; so  $\rho (z,x) \geq 1/(n+1) =  \rho
(x,y)$, and the inequality holds.  If $z$ is not in $U'$, then $z$
and $y$ are in different members of ${\mathcal U}_n$ so as before
the inequality holds.
  The only remaining possibility is that $z  \in    U  \cap  U'$,
but this is impossible because $\calU_n$ is a pairwise disjoint
family.
    To finish the proof we have to show that  $\rho$  is topologically
equivalent to the original metric $d$ on $X$.  It is clear that any
sequence of families of open sets in $(X,d)$ satisfying (i) is a
base for $(X,d)$.  By definition of the metric topology, $\{N_{\rho
}(x,r): x  \in    X, \mbox{ and }r > 0\}$ is a base for $(X, \rho
)$; so to see that $T(X,d) = T(X, \rho )$ it suffices to show that

\begin{equation}\label{equality}\{N_{\rho}(x,r): x  \in
X \mbox{ and }  r > 0\}  =  \cup\{{\mathcal U}_n: n  \in \omega\}
\cup\{X\}.\end{equation}

To this end, let $N = N_{\rho }(x,r)$ be given.  If $r > 1$, then $N
= X$.  If $r \leq 1$, there exists $n \geq 1$ such that $1/(n+1) < r
\leq  1/n$.  Since ${\mathcal U}_{n-1}$ covers $X$, there exists $U
\in {\mathcal U}_{n-1}$ such that $x  \in    U$.  We claim that $U =
N$. If $y  \in    U$, then  $x$ and $y$ are in the same element of
$\calU_{n-1}$ so by (ii), $k(x,y) \geq n$. Hence  $\rho(x,y) \leq
1/(n+1) < r$; so $y  \in    N$.  Conversely, if $ y  \in N$, and $y
\not= x$, then let $m  \in    \omega$ be such that $k(x,y) = m$.
Since $y  \in    N$, we have  $\rho (x,y) = 1/(m+1) < r \leq1/n$; so
$n < m+1$, hence $n-1 < k(x,y)$.  Thus $x$ and $y$ are in the same
member of $\calU_{n-1}$; so $y  \in    U$. To complete the proof of
the above equality (\ref{equality}), let $n  \in    \omega$ and $U
\in    {\mathcal U}_n$.  We claim that $U = N_{\rho}(x,1/(n+1))$ for
any $x$ in $U$.  If $y  \in    U$, then $x$ and $y$ are in the same
element of $U_n$; so $k(x,y) \geq n+1$; so  $\rho (x,y) \leq 1/(n+2)
< 1/(n+1)$, and thus $y  \in    N_{\rho}(x,1/(n+1))$.  If $y  \in
N_{\rho}(x,1/(n+1))$, then $\rho(x,y) \leq 1/(n+2)$; so $k(x,y) \geq
n+1$.  Hence $x$ and $y$ are in the same member of $\calU_n$, and
thus $y \in  U$.  This completes the proof of Theorem
\ref{ultrametricchar}.

\begin{remarks}{\em

Some parts of Theorem \ref{ultrametric} were first noted for the
case of the Baire space in \cite{FellerTornier}

Our proof of Theorem \ref{ultrametricchar} defines the ultrametric
directly on the set $X$ rather than embedding $X$ into a universal
ultrametric space.

The $p$-adic valuations on the rational numbers (for prime numbers
$p$) are examples of other ``naturally occurring" ultrametric spaces
(e.g., see G. Bachman \cite{bachman}).  Indeed the notion of the
strong triangle inequality seems to have arisen from the study of
valuation fields as early as 1918 \cite{ostrowski}. }

\end{remarks}

\section{Continued Fraction Homeomorphism}
In this section we discuss the continued fraction homeomorphism
between Baire space and the irrationals. We use the continued
fractions to define end-points of intervals.  Working by induction,
the  open intervals with these end-points are collected to form a
sequence of open covers $\{\calu_i:i\in\omega\}$ of the irrationals.
We then use these open covers to define the homeomorphism.

From one point of view these open covers are moderately simple to
visualize.  Each $\calU_i$ consists of countably many pairwise
disjoint open intervals with rational end-points. To construct
$\calu_{i+1}$ from $\calU_i$, we take each interval $I = (a,b) \in
\calU_i$, and partition it into countably many intervals $I_j$ ($j
\in \omega$) in the following way:  We construct a sequence of
rational numbers $(p_n)$ that begins with one of the end-points
($p_0 = b$ when $i$ is even; $p_0 = a$ when $i$ is odd) and
converges monotonically to the other end-point.  Consecutive
rational numbers in this sequence are used as the end-points of
intervals in $\calU_{i+1}$ (we define the sequence $(p_n)$ using
continued fractions).

Using continued fractions introduces a bit of  algebra that is
unnecessary for the topology (from the topological point of view, it
does not matter how the sequence $(p_n)$ is defined), but this
construction
 serves as an introduction to  continued fractions,
which is a topic of considerable interest itself.
 Furthermore, a number of theorems in classical topology were originally proved
using continued fractions in much the way presented here.  For
example, the theorem of Sierpi\'nski at the end of this section, and
the well-known topological characterization of the irrationals, due
to Alexandroff and Urysohn \cite{AU}.

\subsection{Continued Fractions}

An expression of the following form is called a (simple) continued
fraction:

\begin{center}$
\def\cfrac#1#2{{\displaystyle\strut#1\over\displaystyle#2}%
 \kern-\nulldelimiterspace}
a_0 + \cfrac{1}{a_1+
 \cfrac{1}{a_2+
  \cfrac{1}{a_3+\cdots_{\vdots\cdots
   \cfrac{1}{x}}
   }}} = [a_0,a_1,\cdots,x].
$
\end{center}
We only consider the special case where $a_0\in \BBZ$,  $a_i\in
\BBN$ for all $1\leq i < n$, and $x$ is any positive number $x \not=
1$ (if $x = 1$, we write $[a_0,a_1,\cdots,a_{n-1}+1]$ instead of
$[a_0,a_1,\cdots,a_{n-1},1]$). The notation $[a_0,a_1,\cdots,x]$  is
commonly used to denote the above continued fraction, where $[a_0] =
a_0$.
  We note the following three
obvious properties: (1)  If the last term $x$ in $[a_0,...,x] $ is a
positive integer $x = a_n$,
 then the number $[a_0,...,a_n] $ is a rational number, (2)
$[a_0,...,a_n] = a_0 + 1/[a_1,\cdots ,a_n]$, (3) for any number
$x\not= 0$,
$$[a_0,\cdots,a_{n-1}, x] = [a_0,\cdots,a_{n-2},a_{n-1} + \frac{1}{x}].$$

\medskip
The continued fraction expansion of a rational number $\frac{a}{b}$
can be found by using the Euclidean algorithm, and is easily seen to
be unique (exercise).

\begin{theorem}[Euclid] For any two integers $a,b$ ($b > 0$) there exists  integers $q$
(the quotient) and $r$ (the remainder) such that
\[a = b\,q + r\mbox{ where } 0 \leq r < b.\]
\end{theorem}

Given two integers $a,b$ we  apply the Euclidean algorithm
repeatedly and set up a sequence of equalities as follows (where we
write $b = r_0$ to facilitate the notation)

\[a = r_0\,q_0 + r_1\mbox{ where } 0 \leq r_1 < r_0\]
\[r_0 = r_1\,q_1 + r_2\mbox{ where } 0 \leq r_2 < r_1\]
\[\vdots\]
\[\label{last}r_{n-2} = r_{n-1}\,q_{n-1} + r_n\mbox{ where }
0 \leq r_n < r_{n-1}.\]
\[\vdots\]

Since the $(r_i)$ form a decreasing sequence of non-negative
integers, this process will end in a finite number of steps. That is
to say, there exists some integer $n$ such that $r_n = 0$.  Our
interest in the above system of equalities is that we can pull the
continued fraction expansion of $ \frac{a}{b}$ from it.  Indeed,

\begin{center}$
\def\cfrac#1#2{{\displaystyle\strut#1\over\displaystyle#2}%
 \kern-\nulldelimiterspace}
\frac{a}{b} = q_0 + \cfrac{1}{q_1+
 \cfrac{1}{q_2+
  \cfrac{1}{q_3+\cdots_{\vdots\cdots
   \cfrac{1}{q_n}}
   }}} = [q_0,q_1,\cdots,q_n]
$
\end{center}

To see this, we rewrite the above family of equalities as follows.

\[\frac{a}{b} = \frac{a}{r_0} = q_0 + \frac{r_1}{r_0} \]
\[\frac{r_0}{r_1} = q_1 + \frac{r_2}{r_1}\mbox{,\,\,\,\,\,or\,\,\,\ }
\frac{r_1}{r_0} = \frac{1}{q_1 + \frac{r_2}{r_1}}\]
\[\vdots\]
\[\label{last1}\frac{r_{n-2}}{r_{n-1}} = q_{n-1}\mbox{,\,\,\,\,\,or\,\,\,\ }
\frac{r_{n-1}}{r_{n-2}}=\frac{1}{q_{n-1}}. \]

\begin{theorem}  Any rational number can be represented as a finite
continued fraction (and conversely, every finite continued fraction
is a rational number)
\end{theorem}

We now give several lemmas concerning continued fractions.

\begin{lemma}\label{x<y}If $x < y$ then the following hold for all $n \geq 0$

(i) If $n$ is even then, $[a_0,a_1,\cdots,a_{n},x] >
[a_0,a_1,\cdots,a_{n},y].$

(ii) If $n$ is odd, $[a_0,a_1,\cdots,a_{n},x] <
[a_0,a_1,\cdots,a_{n},y].$
\end{lemma}

Proof (by induction). For the case $n = 0$ we have $[a_0,x] =
a_0+\frac{1}{x} >  a_0+\frac{1}{y} = [a_0,y]$.  For $n = 1$, we have
$a_1 +\frac{1}{x} > a_1 + \frac{1}{y}$; so by the case $n = 0$ we
have
$$[a_0,a_1,x] = [a_0,a_1+\frac{1}{x}] < [a_0,a_1+\frac{1}{y}] = [a_0,a_1,y]$$
Proof of (i).  Assume the result holds for $2n$ and $2n+1$.  Let
$x^{\prime} = a_{2n+2} + \frac{1}{x}$, and $y^{\prime} = a_{2n+2} +
\frac{1}{y}$; thus $y^{\prime} < x^{\prime}$. Now we have by (ii)
 $$[a_0,a_1,\cdots,a_{2n+2},y] = [a_0,a_1,\cdots,a_{2n+1},y^{\prime}] <
[a_0,a_1,\cdots,a_{2n+1},x^{\prime}] =
[a_0,a_1,\cdots,a_{2n+2},x].$$

Proof of (ii).  Let
 $x^{\prime} = a_{2n+3} + \frac{1}{x}$, and
$y^{\prime} = a_{2n+3} + \frac{1}{y}$; thus $y^{\prime} <
x^{\prime}$.  By (i) just proved, we have
 $$[a_0,a_1,\cdots,a_{2n+3},x] = [a_0,a_1,\cdots,a_{2n+2},x^{\prime}] <
[a_0,a_1,\cdots,a_{2n+2},y^{\prime}] =
[a_0,a_1,\cdots,a_{2n+3},y].$$

\begin{corollary}\label{cor}For any $k\in N$,

(i) if $n$ is even, $[a_0,a_1,\cdots,a_{n},k] >
[a_0,a_1,\cdots,a_{n},k+1]$,

(ii) if $n$ is odd, $[a_0,a_1,\cdots,a_{n},k] <
[a_0,a_1,\cdots,a_{n},k+1]$.
\end{corollary}

\medskip
From this corollary we see that the sequence
$\{[a_0,a_1,\cdots,a_{n},k]: k \in \BBN\}$ is monotone (decreasing
for $n$ even, and increasing for $n$ odd).  We next check that this
sequence converges to $[a_0,a_1,\cdots,a_{n}]$.

\begin{lemma}\label{converges}The sequence $\{[a_0,a_1,\cdots,a_{n},k]: k \in \BBN\}$ converges to
 $[a_0,a_1,\cdots,a_{n}]$.
\end{lemma}

Proof.  Clearly $[a_0,k] = a_0 + \frac{1}{k} \rightarrow a_0 =
[a_0]$, and $[a_0,a_1,k] = a_0 + 1/[a_1,k] \rightarrow a_0+ 1/a_1$
(by case $n = 0$). The result now follows by induction from the
equality $[a_0,...,a_n,k] = a_0 + 1/[a_1,\cdots ,a_n,k].$

\medskip
Now we define a sequence  $\{\calU_i: i < \omega\}$ of  families of
open intervals with rational end-points as discussed in the
introduction to this section.  Put

 $$\calu_0 =\{(a_0,a_0+1): a_0 \in \BBZ\}, \mbox{ and }
\calu_1 =\{([a_0,k+1],[a_0,k]): a_0 \in \BBZ, k\geq 1\}.$$

\noindent Assume that we have defined $\calu_{2n+1}$, and that every
interval in $\calu_{2n+1}$ is of the form
$([a_0,a_1,\cdots,a_{2n},k+1],[a_0,a_1,\cdots,a_{2n},k])$.  For
every $(a_0,a_1,\cdots,a_{2n+1})$ define
$\calu_{2n+2}((a_0,a_1,\cdots,a_{2n+1})) $ to be the following
countable set of intervals:
$$
\{([a_0,a_1,\cdots,a_{2n+1},k],[a_0,a_1,\cdots,a_{2n+1},k+1]): a_0
\in \BBZ, a_i \in \BBN (i \leq 2n+1), k\geq 1\}
$$
By Corollary \ref{cor}, $\{[a_0,a_1,\cdots,a_{2n+1},k]:k\geq 1\}$ is
an increasing sequence with smallest element
$[a_0,a_1,\cdots,a_{2n+1},1] = [a_0,a_1,\cdots,a_{2n+1}+1]$, and by
\ref{converges} converges up to $[a_0,a_1,\cdots,a_{2n+1}]$.  Thus
every interval in $\calu_{2n+2}((a_0,a_1,\cdots,a_{2n+1}))$ is a
subset of the interval
$$I = ([a_0,a_1,\cdots,a_{2n+1}+1], [a_0,a_1,\cdots,a_{2n+1}]) \in  \calu_{2n+1}.$$
Since every member of this increasing sequence (except the smallest
member) is a number strictly between the end-points of $I$, every
interval in $\calu_{2n+2}((a_0,a_1,\cdots,a_{2n+1}))$ has its
closure a subset of $I$ too, except for that one interval which has
$[a_0,a_1,~\cdots,~a_{2n+1}+1~]$ as an end-point.

Define

$$\calu_{2n+2} = \cup\{\calu_{2n+2}((a_0,a_1,\cdots,a_{2n+1}):
 a_0 \in \BBZ,
a_i \in \BBN (i \leq 2n+1)\}$$

In a similar manner, define $\calu_{2n+3}(a_0,a_1,\cdots,a_{2n+2})$
for each  $(a_0,a_1,\cdots,a_{2n+2})$, and define $\calu_{2n+3}$ as
the union over all these countable families of intervals. Since the
sequence
 $ \{[a_0,a_1,\cdots,a_{2n+2},k]:k\geq 1\}$ is decreasing and  converges
down to $[a_0,a_1,\cdots,a_{2n+2}]$ we see that every interval in
$\calu_{2n+3}((a_0,a_1,\cdots,a_{2n+2})$ is a subset of the interval
$$I^{\prime} = ([a_0,a_1,\cdots,a_{2n+2}],[a_0,a_1,\cdots,a_{2n+2}+1])\in \calu_{2n+2}.$$
Moreover, since every member of this decreasing sequence is a member
of the open interval $I^{\prime}$, except the largest member, every
interval in $\calu_{2n+3}((a_0,a_1,\cdots,a_{2n+2}))$ has its
closure a subset of $I^{\prime}$ too, except for that one interval
which has $[a_0,a_1,~\cdots,~a_{2n+2}+1~]$ as an end-point.

In general, since no end-point  of any interval in the  cover
$\calu_{j+2}$ is used as an end-point for any interval in
$\calu_{j}$, we have that every interval in $\calu_{j+2}$ has its
closure a subset of some interval in $\calu_{j}$.

\medskip
We want to show that $\lim_{n \rightarrow
\infty}(\mbox{\,mesh\,}\calu_n) = 0$. We first give a  estimate to
the distance between two continued fractions.  This estimate is
rough, but adequate for our purposes.

\begin{lemma}\label{estimate}Let $1 \leq x < y$ and  $n\geq 1$,

(i) if $n$ is even,
$$[a_0,a_1,\cdots,a_{n},x] -  [a_0,a_1,\cdots,a_{n},y] < \frac{y-x}{xy + n}$$
(ii) if $n$ is odd,
$$[a_0,a_1,\cdots,a_{n},y] -  [a_0,a_1,\cdots,a_{n},x] < \frac{y-x}{xy + n}$$
\end{lemma}
 Proof. For $n = 0$ a quick check shows that $[a_0,x]-[a_0,y] = \frac{y-x}{xy}$.
 We leave the case $n = 1$ as an exercise, and give the inductive step.
Assume the result holds for $2n$ and $2n+1$.  For (i),  put $X =
[a_0,a_1,\cdots,a_{2n+2},x]$ and $Y = [a_0,a_1,\cdots,a_{2n+2},y]$,
and let $x^{\prime} = a_{2n+2}+\frac{1}{x}$ and $y^{\prime} =
a_{2n+2}+\frac{1}{y}$ (thus $y^{\prime}<x^{\prime}$). By (ii) and
the facts that $1 \leq a_{2n+1}, x, y$ we have

\begin{eqnarray*}
X -  Y & = & [a_0,a_1,\cdots,a_{2n+2},x] - [a_0,a_1,\cdots,a_{2n+2},y] \\
& = & [a_0,a_1,\cdots,a_{2n+1},x^{\prime}] -  [a_0,a_1,\cdots,a_{2n+1},y^{\prime}] \\
 & \leq & \frac{x^{\prime}-y^{\prime}}{x^{\prime} y^{\prime}+2n+1}\\
& = & \frac{(a_{2n+1} + \frac{1}{x})-(a_{2n+1} +
\frac{1}{y})}{(a_{2n+1} + \frac{1}{y})
(a_{2n+1}+\frac{1}{x})+2n+1}\\
& = &\frac{y - x}{xya_{2n+1}^2 + a_{2n+1}(x + y) + (2n+1)xy}\\
& < & \frac{y-x}{xy + 2n + 2}
\end{eqnarray*}
The proof of (ii) is similar.

\begin{corollary}\label{mesh}For $n \geq 1$,  $mesh\,(\calu_{n}) < \frac{1}{n+1}$.
\end{corollary}

Proof. The mesh of $\calu_0$ is obviously 1, and the mesh of
$\calu_1$ is easily seen to be $1/2$. Let  $I$ be an arbitrary
interval in $\calu_{n}$ $(n \geq 2)$, $n$ even, say

$$I = ([a_0,a_1,\cdots,a_{n},k+1],[a_0,a_1,\cdots,a_{n},k])$$
thus, the length of $I$ is
$$L = [a_0,a_1,\cdots,a_{n},k]-[a_0,a_1,\cdots,a_{n},k+1]$$

By Lemma \ref{estimate} we get

$$L <\frac{(k+1) - k}{k(k+1) + n}
< \frac{1}{1+n}.$$ The case for $n$ odd is similar.

\medskip
We  summarize the properties of the covers $\calu_i$.

\begin{remark}\label{summary}  The following properties hold for the preceding construction.

(1)  Each $\calu_i$ is a family of open intervals with rational
end-points,  and $\calu_i$ covers $\BBP$,

(2)  The intervals in $\calu_i$ are pairwise disjoint,

(3) $\calu_{i+1}$ refines $\calu_i$, i.e., each interval in
$\calu_{i+1}$ is a subset of (a necessarily unique) interval in
$\calu_i$.

(4)  $\overline{\calu_{i+2}}$ refines $\calu_i$, i.e., the closure
of each interval in $\calu_{i+2}$ is a subset of some interval in
$\calu_i$.

(5)   For $n > 0$, mesh$(\,\calu_i) < \frac{1}{i+1}$.
\end{remark}

Proof. Properties (1) - (4) follow from the construction of the
$\calu_i$, and (5) is Corollary \ref{mesh}.

\begin{corollary}\label{onto}
For every irrational number $x$, there is a sequence
$$\{a_i: i\in \omega,a_0 \in \BBZ
\mbox{, and }a_i\in \BBN\mbox{ for all }i\in \BBN\}$$ such that

1.  (even case) $x$ is the unique point in
\begin{center}$\cap\{([a_0,a_1,\cdots,a_n],[a_0,a_1,\cdots,a_n+1]):n\in\omega, n\mbox{ even}\},
\mbox{ and}$
\end{center}

2. (odd case), $x$ is the unique point in
\begin{center}$\cap\{([a_0,a_1,\cdots,a_n+1],[a_0,a_1,\cdots,a_n]):n\in\omega,n\mbox{ odd}\}.$
\end{center}\end{corollary}

Proof. Let $x$ be an irrational number.  By \ref{summary} (1), there
exists interval $I_0 \in \calu_0$ such that $x \in I_0$.  By the
construction, there exists an integer $a_0$ such that $I_0 =
(a_0,a_0+1)$. By induction, we assume that we have constructed (for
$i \leq n$) intervals $I_i\in\calu_i$, and positive integers $a_i$
(for $1 \leq i < n)$ such that

1.  $x \in I_i \subset I_{i-1}$

2.  $I_i =([a_0,a_1,\cdots,a_i],[a_0,a_1,\cdots,a_{i}+1])$  for $i$
even, and $I_i =([a_0,a_1,\cdots,a_i+1],[a_0,a_1,\cdots,a_{i}])$ for
$i$ odd. To find $I_{n+1}$ is trivial since there exists exactly one
interval in $\calu_{n+1}$ which contains $x$.  Moreover,
$I_{n+1}\subset I_{n}$; so if $n$ is even, by the construction,
$I_{n+1}$ is an interval defined by constructing a sequence that
begins with the left end-point of $I_n$ and converges down to the
right end-point of $I_n$.  Thus there exists a positive integer $k$
such that  $I_{n+1}=
([a_0,a_1,\cdots,a_n,k+1],[a_0,a_1,\cdots,a_n,k])$. Put $a_{n+1} =
k$. The case that $n$ is odd is similar.

\subsection{The homeomorphism}

We want to construct a homeomorphism $\phi$ from $\BBB$ onto $\BBP$.
We work, however, with $\BBB_2$ which is homeomorphic to $\BBB$ by
Theorem \ref{prod}.

\begin{definition}$\phi((a_i))$ is the unique point in
$$\cap\{([a_0,a_1,\cdots,a_n],[a_0,a_1,\cdots,a_n+1]):n\mbox{ even}\}.$$
\end{definition}

Obviously, we also have
$$\phi((a_i)) \mbox{ is the unique point in }
\cap\{([a_0,a_1,\cdots,a_n+1],[a_0,a_1,\cdots,a_n]):n\mbox{
odd}\}.$$ Let us check that $\phi$ is well-defined.  The
intersection that defines $\phi$ is non-empty because the
intersection of these open intervals equals the intersection of
their closures by  \ref{summary} (4) (and therefore the intersection
is non-empty by the local compactness of $\BBR$).  The intersection
is a singleton by \ref{summary} (5). The unique point $x$ in this
intersection is irrational because every rational number has a
continued fraction expansion, and therefore appears as an end-point
of an interval in some $\calu_i$.  Thus by \ref{summary} (4), $x$ is
not  such an end-point and therefore is irrational.

\begin{theorem}[Baire \cite{baire}] $\phi$ is a homeomorphism from
$\BBB_2$ onto $\BBP$.\end{theorem}

Proof. $\phi$ is one-one:  If $(a_i)\not= (b_i)$, then there exists
a first $j$ such that $a_j \not= b_j$.  If $j = 0$ then $a_0 \not=
b_0$ and the intervals $(a_0,a_0+1)$, and $(b_0,b_0+1)$ are
disjoint; so $\phi((a_i))\not= \phi((b_i))$. In case $j >0$ and even
(so $j-1$ is odd) we have

\[\phi((a_i))\in([a_0,a_1,\cdots,a_{j-1},a_j],[a_0,a_1,\cdots,a_{j-1},a_j+1])\]and
\[\phi((b_i))\in([a_0,a_1,\cdots,a_{j-1},b_j],[a_0,a_1,\cdots,a_{j-1},b_j+1]).\]
Since $a_j\not= b_j$, the two intervals above are distinct elements
of $\calu_j$, hence disjoint.  It follows that $\phi((a_i))\not=
\phi((b_i))$. The case for $j>0$ and odd is similar.

$\phi$ is onto:  This follows at once from \ref{onto}.

$\phi$ is bicontinuous:  It suffices to show that there exist a base
${\mathcal B}$ for the Baire space, and a base  ${\mathcal U}$ for
the irrational numbers such that the mapping $\phi$ sends each $B\in
{\mathcal B}$ onto some $U \in {\mathcal U}$, and each $U\in
{\mathcal U}$ is the image of some $B \in {\mathcal B}$.  For
${\mathcal B}$ we take a usual base for the metric topology on
$\BBB_2$,

$${\mathcal N} = \{N(x,1/n): n \in \BBN\},$$
and for the base ${\mathcal U}$ on the irrationals we take
${\mathcal U} = \cup\{{\mathcal U}_i: i < \omega\}$, which is
clearly a base by \ref{summary} (1), (5).  When using  ${\mathcal
N}$ in $\BBB_2$, it is useful to note that
$$d_2(x,y) < \frac{1}{n} \mbox{ if and only if } x_i = y_i \mbox{ for }i < n.$$
Now it suffices to show that for every $a \in \BBB_2$ and $B =
N(a,1/n) = \{f\in \BBB_2:f_i = a_i \mbox{ for } i < n\}$ we have

\smallskip
(i)  for $n$ even, $\phi(B) =
([a_0,a_1,\cdots,a_n],[a_0,a_1,\cdots,a_n +1])$.

\smallskip
(ii) for $n$ odd, $\phi(B) =
([a_0,a_1,\cdots,a_n+1],[a_0,a_1,\cdots,a_n])$.

\smallskip
\noindent Proof of (i).  For $n$ even, by the definition of $\phi$
and the nested property of the intervals, we have
\begin{eqnarray*}x\in \phi(B)
& \mbox{ iff }& \exists b \in B, b_i = a_i\, (i \leq n) \mbox{ and }\phi(b) = x\\
& \mbox{ iff } &
x \in \cap\{([a_0,a_1,\cdots,a_i],[a_0,a_1,\cdots,a_{i-1},a_i+1]):i\mbox{ even} \leq n\}\\
& \mbox{ iff } & x \in
([a_0,a_1,\cdots,a_n],[a_0,a_1,\cdots,a_{n-1},a_n+1]).
\end{eqnarray*}
The proof of (ii) is similar.
\medskip

One of the many useful features of the Baire space is that for some
purposes it can be used in lieu of the space of irrationals. The
following theorem illustrates this feature.  The theorem is a
special case of a  general theorem of Wac{\l}aw Sierpi\'nski
\cite[Theorem 2]{sierpinski} in which Sierpi\'nski  used  continued
fractions to define a homeomorphism into the irrational numbers.

\begin{theorem}[Sierpi\'nski] Every separable, 0-dimensional  metric space
is homeomorphic to a subset of the irrational numbers.
\end{theorem}

Proof.  Let $X$ be a separable, 0-dimensional  metric space. It
suffices to prove that $X$ is homeomorphic to a subset of $\BBB$. By
Theorem \ref{0covering},we can find a family
$\{\calu_i:i\in\omega\}$ of clopen covers of $X$ satisfying the
conditions in Theorem \ref{ultrametricchar} (3).  By separability,
each $\calU_i$ is at most countable; say has cardinality $k_i$,
where $k_i \leq \omega$ for all $i \in \omega$.  Label the elements
of $\calu_i$ as $\{U_{i,j}: j < k_i\}$ in a one-one manner. For
every $x\in X$ and every $i$ there is a unique
 $U\in\calu_i$ so that $x \in U$,
hence we may pick a unique integer $f_x(i)$ so that $x \in
U_{i,f_x(i)}$. This defines a function $f_x$ where
$$f_x\in\prod_{i<\omega}k_i =\{f\in \BBB: f(i) < k_i\mbox{ for all }i < \omega\},$$
and
$$x \in \cap\{U_{i,f_x(i)}: i \in \omega\}.$$
We claim that the function $\phi:X\rightarrow \BBB$ defined by
$\phi(x) = f_x$ is a homeomorphism onto the image of $\phi$.

$\phi$ is one-one: if $x \not= y$ then there exists $i$ such that
$x,y$ are in different elements of $\calu_i$; so $f_x(i)\not=
f_y(i)$.

$\phi$ is a homeomorphism:  First note that for $i,j\in \omega$,
$\sigma(i,j) = \{f\in\BBB: f(i) = j\}$, is an open set in $\BBB$,
and moreover, the set of all  $\sigma(i,j)$ for $i,j\in\omega$ is a
subbase  for in the  metric topology on $\BBB$. Therefore the family
$$\{\sigma(i,j)\cap\phi(X): i \in \omega, j < k_i\}$$
forms a subbase for $\phi(X)$ since $\phi(x) = f_x \in \prod_{i <
\omega}k_i$.  Now it suffices to show that every set in the base
$\cup\{\calu_i:i < \omega\}$ for $X$ is mapped by $\phi$ onto such a
subbasic set, and every subbasic set is the image of some $U \in
\cup\{\calu_i:i < \omega\}$.  Thus it suffices to prove
$$\phi(U_{i,j}) = [\langle i,j\rangle]\cap\phi(X),$$
for $i < \omega$, and $j < k_i$. This follows immediately from the
definitions involved and is left as an exercise.


\begin{thebibliography}{99}

\bibitem{AU}Paul Alexandroff and Paul Urysohn, \"Uber nulldimensionale Punktmegen,
Math Annalen (1928) 89-106.

\bibitem{bachman}G. Bachman, {\bf Introduction to p-adic numbers and valuation theory},
Academic Press, New York, 1964.

\bibitem{baire}R. Baire, Sur la repr\'esentation des fonctions discontinues
(deuxi\'eme partie), Acta Math. 32 (1909) 97-176.

\bibitem{FellerTornier} Willy Feller and Erhard Tornier, Mas- und Inhaltstheorie
des Baireschen, Math Annalen (1933) 165-187. Nullraumes

\bibitem{kulesza} John Kulesza, An example in the dimension theory of metric spaces,
Topology and Appl. 35 (1990) no. 2-3, 108-120.

\bibitem{adam} Adam Ostaszewski, A note on Prabir Roy's space, Topology and Appl.
35 (1990) no. 2-3, 85-107.

\bibitem{ostrowski}A. Ostrowski, \"Uber einige l\"osungen der funktionalgleichung, Acta Math 41 (1918) 271 -284.

\bibitem{roy}Prabir Roy, Non equality of dimension for metric spaces, Trans. Amer. Math. Soc.
134 (1968) 117-132.


\bibitem{sierpinski}Wac{\l}aw Sierpi\'nski, Sur les ensembles connexes
et non connexes, Fund. Math. 2 (1921) 81 - 95.



\end{thebibliography}
\end{document}